\documentclass[12pt,oneside,reqno]{amsart}
\usepackage{amssymb}
\newtheorem{thm}{Theorem}
\newtheorem{lem}[thm]{Lemma}
\newtheorem{pro}[thm]{Proposition}
\newtheorem{col}[thm]{Corollary}

\newenvironment{remark}[1][Remarks]{\begin{trivlist}
   \item[\hskip \labelsep {\bfseries #1}]}{\end{trivlist}}

\begin{document}
\title[Lattice action on the boundary]{Lattice action on the boundary of $\hbox{\rm SL}(n,\mathbb{R})$}
\author{Alexander Gorodnik}
\address{Department of Mathematics, University of Michigan, Ann Arbor, MI 48109}
\email{gorodnik@umich.edu}
\thanks{This article is a part of author's PhD thesis at Ohio State University done under supervision of Prof.~Bergelson.}

\begin{abstract}
Let $\Gamma$ be a lattice in $G=\hbox{\rm SL}(n,\mathbb{R})$ and $X=G/S$ a homogeneous space of $G$,
where $S$ is a closed subgroup of $G$ which contains a real algebraic subgroup $H$ such that
$G/H$ is compact.  We establish uniform distribution of orbits of $\Gamma$ in $X$ analogous
to the classical equidistribution on torus. To obtain this result, we first prove
an ergodic theorem along balls in the connected component of Borel subgroup of $G$.
\end{abstract}

\maketitle

\section{Introduction}

Let $G=\hbox{\rm SL}(n,\mathbb{R})$ and $\Gamma$ a lattice in $G$; that is, $\Gamma$ is a discrete subgroup of $G$ with finite covolume.
We study the action of $\Gamma$ on a compact homogeneous space $X$ of algebraic origin.
Namely, $X=G/S$ where $S$ is a closed subgroup of $G$ which contains the connected
component of a real algebraic subgroup $H$ of $G$ such that $G/H$ is compact.
An important example is provided by the Furstenberg boundary of $G$ \cite{f63}. In this case, $X=G/B$
where $B$ is the subgroup of upper triangular matrices in $G$.

It is possible to deduce from a result of Dani \cite[Theorem~13.1]{st}
that every orbit of $\Gamma$ in $X$ is dense.
We will prove a quantitative estimate for the distribution of orbits.

Introduce a norm on $G$:
\begin{equation}\label{eq_norm0}
\|g\|=\left(\hbox{Tr}({}^tgg)\right)^{1/2}=\left(\sum_{i,j} g_{ij}^2\right)^{1/2}\quad\hbox{for}\quad g=(g_{ij})\in G.
\end{equation}
For $T>0$, $\Omega\subseteq X$, and $x_0\in X$, define a counting function
\begin{equation}\label{eq_NT}
N_T(\Omega, x_0)=|\{\gamma\in \Gamma: \|\gamma\|<T, \gamma\cdot x_0\in \Omega\}|.
\end{equation}
Let $m$ be a normalized $\hbox{SO}(n)$-invariant measure on $X$.
It follows from the Iwasawa decomposition (see (\ref{eq_iwasawa}) below) that
$X$ is a homogeneous space of $\hbox{SO}(n)$.
Therefore, the measure $m$ is unique.
The following  theorem shows that orbits of $\Gamma$ in $X$ are uniformly distributed with respect
to the measure $m$.

\begin{thm} \label{th_asympt}
For a Borel set $\Omega\subseteq X$ such that $m(\partial \Omega)=0$ and $x_0\in X$,
\begin{equation}\label{eq_asympt}
N_T(\Omega, x_0)\sim \frac{\gamma_n}{\bar\mu(\Gamma\backslash G)} m(\Omega) T^{n(n-1)}\quad\hbox{as}\quad T\rightarrow\infty,
\end{equation}
where $\gamma_n$ is a constant (defined in (\ref{eq_gamma}) below), and $\bar\mu$ 
is a finite $G$-invariant measure on $\Gamma\backslash G$ (defined in (\ref{eq_mubar})).
\end{thm}

It would be interesting to obtain an estimate for the error term in (\ref{eq_asympt}).
This, however, would demand introducing different techniques than those employed here.

Theorem \ref{th_asympt} is analogous to the result of Ledrappier \cite{l} (see also \cite{n})
who investigated the distribution of dense orbits of a lattice $\Gamma\subset\hbox{\rm SL}(2,\mathbb{R})$
acting on $\mathbb{R}^2$. Ledrappier used the equidistribution property of the horocycle flow.
Similarly, we deduce Theorem \ref{th_asympt} from an equidistribution property of orbits of Borel subgroup.

Denote by $B^o$ the connected component of the upper triangular subgroup of $G$.
To prove Theorem \ref{th_asympt}, we use the following ergodic theorem for the right action of $B^o$
on $\Gamma\backslash G$.

\begin{thm} \label{th_ergodic}
Let $\varrho$ be a right Haar measure on $B^o$, and $\nu$ the normalized $G$-invariant measure
on $\Gamma\backslash G$. Then for any $\tilde{f}\in C_c(\Gamma\backslash G)$ and $y\in \Gamma\backslash G$,
$$
\frac{1}{\varrho(B^o_T)}\int_{B^o_T} \tilde{f}(yb^{-1})d\varrho(b)\rightarrow \int_{\Gamma\backslash G} \tilde{f}d\nu\quad\hbox{as}\quad T\rightarrow\infty,
$$
where $B^o_T=\{b\in B^o:\|b\|<T\}$.
\end{thm}

\begin{remark}[Remarks.]
\hspace{1cm}

\begin{enumerate}
\item[1.] One can consider the analogous limit for a left Haar measure on $B^o$.
In this case, it may happen that the limit is $0$ for some $y\in \Gamma\backslash G$
and all $\tilde{f}\in C_c(\Gamma\backslash G)$ (see Proposition \ref{l_contr}).

\item[2.] Since $B^o$ is solvable (hence, amenable),
one might expect that convergence for a.e. $y\in\Gamma\backslash G$
follows from known ergodic theorems for amenable group actions.
Moreover, since $\nu$ is the only normalized $B^o$-invariant measure on $\Gamma\backslash G$,
one could expect that convergence holds for every $y$.
However, this approach does not work 
because the sets $B^o_{T}$ do not form a F\o lner sequence, and the space $\Gamma\backslash G$
is not compact in general.

\item[3.] To prove Theorem \ref{th_ergodic}, we use Ratner's classification
of ergodic measures for unipotent flows \cite{r1}. In fact, we don't need the full strength of 
this result. Since the subgroup $U$ (defined in (\ref{eq_U}) below) is horospherical,
it is enough to know classification of ergodic measures for horospherical subgroups.
The situation is much easier in this special case (see \cite[\S 13]{st}).

\item[4.] We expect that analogs of Theorems \ref{th_asympt} and \ref{th_ergodic} hold
for a noncompact semisimple Lie group and its irreducible lattice with balls $B_T^o$
defined by the Riemann metric. The main difficulty here is to show that
the measure of $B_T^o$ is ``concentrated'' on the ``cone'' $B_T^C$ (cf. Lemma \ref{lem_BTC2}). 

\item[5.] It was pointed out by P.~Sarnak that it might be possible to prove the results of this article using
harmonic analysis on $\Gamma\backslash G$. In particular, Corollary \ref{th_distr_sl2} below can be deduced from
the result of Good (Corollary on page 119 of \cite{g}). Note that his method gives an estimate on the error term.
\end{enumerate}
\end{remark}

The paper is arranged as follows. In the next section, we give examples of applications of Theorem \ref{th_asympt}.
In Section \ref{sec_ba} we set up notations and prove some basic lemmas.
Theorem \ref{th_asympt} is deduced from Theorem \ref{th_ergodic} in Section \ref{sec_th1}.
In Sections \ref{sec_uni} and \ref{sec_rep}, we review results on the structure of unipotent flows and prove auxiliary facts about
finite-dimensional representations of $\hbox{\rm SL}(n,\mathbb{R})$. Finally, Theorem \ref{th_ergodic}
is proved in Section \ref{sec_th2}.

\section{Examples} \label{s_ex}

\begin{enumerate}
\item[1.] Let $X=\mathbb{R}\cup\{\infty\}$, which is considered as the boundary of the hyperbolic upper half plane.
The group $G=\hbox{\rm SL}(2,\mathbb{R})$ acts on $X$ by fractional linear transformations:
\begin{equation}\label{eq_frac}
g\cdot x=\frac{ax+b}{cx+d}\quad \hbox{for}\;x\in X,\; g=\left(
\begin{tabular}{cc}
$a$ & $b$\\
$c$ & $d$
\end{tabular}
\right)\in G.
\end{equation}
Let $\Gamma$ be a lattice in $\hbox{\rm SL}(2,\mathbb{R})$.
For $\Omega\subseteq X$ and $x_0\in X$, define the counting function $N_T(\Omega,x_0)$ as in (\ref{eq_NT}).
Its asymptotics can be derived from Theorem \ref{th_asympt}.
Note that the asymptotics of $N_T(X,x_0)$ as $T\rightarrow\infty$ provides a solution of
the so-called hyperbolic circle problem (cf. \cite[p.~266]{ter1} and references therein).

\begin{col} \label{th_distr_sl2}
{\sc (of Theorem \ref{th_asympt})}
For $x_0\in X$ and $-\infty\le a<b\le +\infty$, 
$$
N_T((a,b),x_0)\sim c_{\Gamma}\left(\int_a^b\frac{dt}{1+t^2}\right) T^2\quad\hbox{as}\quad T\rightarrow\infty,
$$
where $c_\Gamma=\frac{1}{2\bar\mu(\Gamma\backslash G)}$ ($\bar\mu$ is the
$G$-invariant measure on $\Gamma\backslash G$ defined in (\ref{eq_mubar})).
\end{col}

\begin{proof}[Proof.]
It is easy to see from (\ref{eq_frac}) that $G$ acts
transitively on $X$, and the 
stabilizer of $\infty$ in $G$ is the group of upper triangular matrices $B$.
Thus, Theorem \ref{th_asympt} is applicable to the space $X$.

Note that 
$$
K=\hbox{SO}(2)=\left\{k_\theta=\left(\begin{tabular}{rr}
$\cos 2\pi\theta$ & $\sin 2\pi\theta$\\
$-\sin 2\pi\theta$ & $\cos 2\pi\theta$
\end{tabular}
\right):\theta\in [0,1)\right\},
$$
and the normalized Haar measure on $K$ is given by $dk=d\theta$.
The measure $m$ on $X$ can be defined as the image of $dk$ under the map
$K\rightarrow X:k\rightarrow k\cdot \infty$. 
By (\ref{eq_frac}), $k_\theta\cdot\infty=-\hbox{ctan}\; 2\pi\theta$.
Then
$$
m((a,b))=\int_{k\cdot\infty\in (a,b)}dk=\mathop{\int_{-\hbox{\small ctan}\; 2\pi\theta\in (a,b)}}_{\theta\in [0,1)} d\theta
=\frac{1}{\pi}\int_a^b\frac{dt}{1+t^2}.
$$
We have used the substitution $t=-\hbox{ctan}\; 2\pi\theta$.

Finally, by Theorem \ref{th_asympt},
$$
N_T((a,b),x_0)\mathop{\sim}_{T\rightarrow\infty} \frac{\gamma_2 m((a,b))}{\bar\mu(\Gamma\backslash G)}T^2
=c_{\Gamma}\left(\int_a^b\frac{dt}{1+t^2}\right) T^2.
$$
Note that $\gamma_2=\frac{\pi}{2}$ by (\ref{eq_gamma}) below.
\end{proof}\medbreak

\item[2.] Let $X=\mathbb{P}^{n-1}$ be the projective space
(more generally $X=\mathcal{G}_{n,k}$, Grassmann variety, or $X=\mathcal{F}_n$,
flag variety), and $m$ be the rotation invariant normalized measure on $X$.
Then the asymptotic estimate (\ref{eq_asympt}) holds for the standard action of $G=\hbox{\rm SL}(n,\mathbb{R})$
on $X$. This is a special case of Theorem \ref{th_asympt} because $X$ can be identified
with $G/S$ where $S$ is a closed subgroup of $G$ that contains $B$, the group of
upper triangular matrices.
\end{enumerate}

\section{Basic facts} \label{sec_ba}

For $s=(s_1,\ldots, s_n)\in\mathbb{R}^n$, define
$$
a(s)=\hbox{diag}(e^{s_1},\ldots, e^{s_n}).
$$
For $t=(t_{ij}:1\le i<j\le n)$, $t_{ij}\in\mathbb{R}$, denote by $n(t)$ the unipotent upper triangular matrix with
entries $t_{ij}$ above the main diagonal.

We use the following notations:
\begin{eqnarray}
G&=&\hbox{\rm SL}(n,\mathbb{R}),\nonumber\\
K&=&\hbox{SO}(n),\nonumber\\
A^o&=&\{a(s)|\;s\in\mathbb{R}^n,\;\sum_i s_i=0\},\nonumber\\
N&=&\{n(t)|\; t_{ij}\in\mathbb{R},1\le i<j\le n\},\nonumber\\
B^o&=&A^oN=NA^o. \nonumber
\end{eqnarray}
For $s\in\mathbb{R}^n$, $\sum_i s_i=0$, denote $\alpha_{ij}(s)=s_i-s_j$, where $i,j=1,\ldots, n$,
$i<j$. These are the positive roots of the Lie algebra of $G$. Note that
\begin{equation} \label{eq_adj}
\hbox{Ad}_{a(s)}n\left(\{t_{ij}\}\right)=a(s)n\left(\{t_{ij}\}\right)a(s)^{-1}=n(\{e^{\alpha_{ij}(s)}t_{ij}\}).
\end{equation}
Let 
\begin{equation} \label{eq_varrho}
\delta(s)=\frac{1}{2}\sum_{i<j} \alpha_{ij}(s)=\sum_{1\le k\le n} (n-k)s_k.
\end{equation}


For $C\in \mathbb{R}$, define
$$
A^C=\{a(s)\in A^o|\;s_i>C,i=1,\ldots, n-1\}.
$$
Also put $B^C=A^CN$.

Let $dk$ be the normalized Haar measure on $K$.
A Haar measure on $N$ is given by $dn=dt=\prod_{i<j}dt_{ij}$.
A Haar measure on $A^o$ is $da=ds=ds_1\ldots ds_{n-1}$.

The product map $A^o\times N\rightarrow B^o$ is a diffeomorphism. 
The image of the product measure under this map is a left Haar measure on $B^o$.
Denote this measure by $\lambda$. Then a right Haar measure $\varrho$ on $B^o$ can be defined by
\begin{equation} \label{eq_lambda}
\varrho (f)=\int_{B^o} f(b^{-1}) \lambda(b)=\int_{A^o\times N} f(a(s)n(t))e^{2\delta(s)}dsdt,\quad f\in C_c(B^o).
\end{equation}

The map corresponding to the Iwasawa decomposition 
\begin{equation} \label{eq_iwasawa}
(k,a,n)\mapsto kan: K\times A^o\times N\rightarrow G
\end{equation}
is a diffeomorphism. One can define a Haar measure $\mu$ on $G$
in terms of this decomposition:
\begin{equation} \label{eq_iwasawa2}
\int_G fd\mu=\int_{K\times A^o\times N}f(ka(s)n(t))e^{2\delta(s)}dkdsdt=\int_{K\times B^o} f(kb)dkd\varrho(b)\end{equation}
for $f\in C_c(G)$.
For a lattice $\Gamma$ in $G$, there exists
a finite measure $\bar \mu_\Gamma$ on $\Gamma\backslash G$ such that
\begin{equation} \label{eq_mubar}
\int_G fd\mu=\int_{\Gamma\backslash G}\sum_{\gamma\in\Gamma} f(\gamma g)d\bar\mu_\Gamma(g),\quad f\in C_c(G).
\end{equation}

Let $\beta$ be an automorphism of $G$. Then $\beta(\Gamma)$ is a lattice too.
Moreover, the following lemma holds.

\begin{lem}\label{l_auto}
Define a map 
$$
\bar\beta:\Gamma\backslash G\rightarrow\beta(\Gamma)\backslash G:g\Gamma\mapsto \beta(g)\beta(\Gamma)
$$
Then $\bar\beta(\bar\mu_\Gamma)=\bar\mu_{\beta(\Gamma)}$.
\end{lem}

\begin{proof}[Proof.]
Since the automorphism group of $G$ is a finite extension of the group of inner automorphisms,
and $G$ is unimodular, it follows that that the measure $\mu$ is $\beta$-invariant.

Every $\tilde{f}\in C_c(\beta(\Gamma)\backslash G)$ can be represented as
$\sum_{\gamma\in\beta(\Gamma)} f(\gamma g)$ for some $f\in C_c(G)$ (see \cite[Ch.~1]{rag}).
Then
\begin{eqnarray*}
\bar\beta(\bar\mu_\Gamma)(\tilde{f})&=&
\int_{\Gamma\backslash G}\sum_{\gamma\in\beta(\Gamma)} f(\gamma \beta(g))d\bar\mu_\Gamma(g)
=\int_G f(\beta(g))d\mu(g)\\
&=&\int_G f(g)d\mu(g)=\bar\mu_\Gamma(\tilde{f}).
\end{eqnarray*}
\end{proof}\medbreak

For a subset $D\subseteq G$ and $T>0$, put
$$
D_T=\{d\in D: \|d\|<T\}.
$$
Note that
\begin{equation} \label{eq_BT}
B_T^o=\left\{a(s)n(t): \sum_{1\le i\le n} e^{2s_i}+\sum_{1\le i<j\le n} e^{2s_i}t_{ij}^2<T^2 \right\}.
\end{equation}
For $s\in\mathbb{R}^n$, define
\begin{equation} \label{eq_Ns}
N(s)=\sum_i e^{2s_i}.
\end{equation}
\begin{lem} \label{lem_BTC}
For $C\in\mathbb{R}$,
$$
\varrho(B_T^C)=c_n\int_{A^C_T} \Big(T^2-N(s)\Big)^{\frac{n(n-1)}{4}}\hbox{\rm exp}\left(\sum_{k} (n-k)s_k\right)ds,
$$
where $c_n=\pi^{n(n-1)/4}/\Gamma(1+n(n-1)/4)$.
\end{lem}

\begin{proof}[Proof.]
Use formulas (\ref{eq_lambda}), (\ref{eq_varrho}), (\ref{eq_BT}), and make change of
variables $t_{ij}\rightarrow e^{-s_i}t_{ij}$. Then the above formula follows from the fact
that the volume of the unit ball in $\mathbb{R}^m$ is $\pi^{m/2}/\Gamma(1+m/2)$.
\end{proof}\medbreak

It follows from Lemma \ref{lem_BTC} that $\varrho (B_T^o)=O\left(T^{(n^2-n)}\right)$ as $T\rightarrow\infty$.
In fact, more precise statement is true:
\begin{lem}
\begin{equation} \label{eq_BTasy}
\varrho(B_T^o)\sim \gamma_n T^{(n^2-n)}\quad\hbox{as}\quad T\rightarrow\infty,
\end{equation}
where
\begin{equation} \label{eq_gamma}
\gamma_n=\frac{\pi^{n(n-1)/4}}{2^{n-1}\Gamma\left(\frac{n^2-n+2}{2}\right)}\prod_{k=1}^{n-1} \Gamma\left(\frac{n-k}{2}\right).
\end{equation}
\end{lem}

\begin{proof}[Proof.]
Since the norm (\ref{eq_norm0}) is $K$-invariant, $G_T=KB^o_T$. By (\ref{eq_iwasawa2}), $\mu (G_T)=\varrho(B^o_T)$.
The asymptotics of $\mu (G_T)$ as $T\rightarrow\infty$ was computed in \cite[Appendix~1]{drs}.
\end{proof}\medbreak

\begin{lem} \label{lem_BTC2}
For any $C\in\mathbb{R}$, $\varrho(B_T^C)\sim\varrho(B_T^o)$ as $T\rightarrow\infty$.
\end{lem}

\begin{proof}[Proof.]
For $i_0=1,\ldots, n-1$, put 
$$
A_T^{i_0}=\{a(s)\in A_T^o: s_{i_0}\le C\}\;\;\hbox{and}\;\;B_T^{i_0}=\{a(s)n(t)\in B_T^o: s_{i_0}\le C\}.
$$
We claim that $\varrho (B_T^{i_0})=o(\varrho (B_T^o))$ as $T\rightarrow\infty$.
It follows from (\ref{eq_BT}) that if $a(s)\in A_T^o$, then
$s_i<\log T$ for every $i=1,\ldots, n$.
Then by Lemma \ref{lem_BTC}, 
\begin{eqnarray} 
\nonumber \varrho(B_T^{i_0})&\le& c_nT^{\frac{n(n-1)}{2}}\int_{A_T^{i_0}} \hbox{\rm exp}\left(\sum_{k} (n-k)s_k\right)ds \\
\nonumber &\ll& T^{\frac{n(n-1)}{2}} \prod_{{k<n},{k\ne i_0}} \int_{-\infty}^{\log T} e^{(n-k)s_k}ds_k\ll T^{n(n-1)-(n-i_0)}.
\end{eqnarray}
(Here and later on $A\ll B$ means $A < c\cdot B$ for some absolute constant $c>0$.)
Now the claim follows from (\ref{eq_BTasy}).
Since $B_T^o-B_T^C=\cup_{i_0<n} B_T^{i_0}$, $\varrho (B_T^o-B_T^C)=o(\varrho (B_T^o))$
as $T\rightarrow\infty$. Therefore, $\varrho(B_T^C)\sim\varrho(B_T^o)$ as $T\rightarrow\infty$.
\end{proof}\medbreak

Next, we show that Theorem \ref{th_ergodic} fails for a left Haar measure on $B^o$.

\begin{pro}\label{l_contr}
Let $\Gamma$ be a lattice in $G=\hbox{\rm SL}(2,\mathbb{R})$, and $y\in\Gamma\backslash G$ be such that
the orbit $yN$ is periodic. Then for every $\tilde{f}\in C_c(\Gamma\backslash G)$,
$$
\frac{1}{\lambda(B^o_T)}\int_{B^o_T} \tilde{f}(yb^{-1})d\lambda(b)\rightarrow 0 \quad\hbox{as}\quad T\rightarrow\infty.
$$
\end{pro}

\begin{proof}[Proof.]
For $C\in \mathbb{R}$, put $\hat{A}^C=\{a(s)\in A^o:s_1<C\}$ and $\hat{B}^C=\hat{A}^C N$.
As in Lemma \ref{lem_BTC2}, one can show that $\lambda (B^o_T-\hat{B}^C_T)=o(\lambda (B^o_T))$
as $T\rightarrow\infty$. Therefore, for every $C\in\mathbb{R}$,
$$
\frac{1}{\lambda(B^o_T)}\int_{B^o_T-\hat{B}^C_T} \tilde{f}(yb^{-1})d\lambda(b)\rightarrow 0 \quad\hbox{as}\quad T\rightarrow\infty.
$$
On the other hand, according to \cite[Lemma~14.2]{st}, $yN^{-1}a(s)^{-1}\rightarrow\infty$ as
$s_1\rightarrow -\infty$. Thus, there exists $C\in\mathbb{R}$ such that 
$y(\hat{B}^C_T)^{-1}\cap \hbox{supp} (\tilde{f})=\emptyset$. Then
$$
\int_{\hat{B}^C_T} \tilde{f}(yb^{-1})d\lambda(b)= 0.
$$
This proves the proposition.
\end{proof}\medbreak

\section{Proof of Theorem \ref{th_asympt}} \label{sec_th1}

{\bf Claim}. {\it It is enough to prove Theorem \ref{th_asympt} for $X=G/B^o$}.

\begin{proof}[Proof.]
Suppose that the theorem is proved for $X=G/B^o$.

At first, we consider a special case: 
\begin{equation}\label{eq_spec}
X=G/(B^o)^{g_0}\quad\hbox{for}\;\;\hbox{some}\;\; g_0\in G,
\end{equation}
where $(B^o)^{g_0}=g_0^{-1}B^og_0$.
By the Iwasawa decomposition (\ref{eq_iwasawa}), $(B^o)^{g_0}=(B^o)^{k_0}$ for some
$k_0\in K$.

The normalized $K$-invariant measure $m$ on $G/(B^o)^{k_0}$ is defined as 
$$
m(C)=\int_{k: k(B^o)^{k_0}\in C} dk\quad\hbox{for}\;\; \hbox{Borel}\;\hbox{set}\;\; C\subseteq G/(B^o)^{k_0}.
$$
Similarly, one defines the normalized $K$-invariant measure $m^*$ on $G/B^o$.
Consider a map
$$
\beta:G/B^o\rightarrow G/(B^o)^{k_0}: gB^o\mapsto g^{k_0} (B^o)^{k_0}.
$$
Clearly, $\beta$ is a diffeomorphism.
Using that $K$ is unimodular, one proves that
\begin{equation}\label{eq_m}
m^*(\beta^{-1}(C))=m(C)\quad\hbox{for}\;\; \hbox{Borel}\;\hbox{set}\;\; C\subseteq G/(B^o)^{k_0}.
\end{equation}
Take 
$$
\Omega\subseteq G/(B^o)^{k_0}\quad\hbox{and}\quad x_0=h_0(B^o)^{k_0}\in G/(B^o)^{k_0}
$$
such that $m(\partial\Omega)=0$. Set 
$$
{\Omega}^*=\beta^{-1}(\Omega)\subseteq G/B^o\quad\hbox{and}\quad {x}^*_0=h_0^{k_0^{-1}}B^o\in G/B^o.
$$
By (\ref{eq_m}), $m^*(\partial{\Omega^*})=0$ too. For $\gamma\in\Gamma$,
$\gamma\cdot x_0\in\Omega$ iff $\gamma^{k_0^{-1}}\cdot {x}^*_0\in {\Omega}^*$.
Note also that $\|\gamma^{k_0^{-1}}\|=\|\gamma\|$. Therefore,
$$
N_T(\Omega,x_0)=|\{\gamma\in \Gamma^{k_0^{-1}}: \|\gamma\|<T, \gamma\cdot {x}^*_0\in {\Omega}^*\}|
$$
Applying the assumption to the lattice $\Gamma^*=\Gamma^{k_0^{-1}}$, one gets
$$
N_T(\Omega, x_0)\sim \frac{\gamma_n}{\bar\mu^*(\Gamma^*\backslash G)} m^*(\Omega^*) T^{n(n-1)}\quad\hbox{as}\quad T\rightarrow\infty,
$$
where $\bar\mu^*$ is the measure on $\Gamma^*\backslash G$ defined in (\ref{eq_mubar}).
Now the special case (\ref{eq_spec}) follows
from Lemma \ref{l_auto} and (\ref{eq_m}).

Let us consider the general case.
Let $S$ be a closed subgroup of $G$ such that $S\supseteq H^o$, where $H$ is a real algebraic subgroup of
$G$, and $G/H$ is compact. Since $H$ has finitely many connected components, $G/H^o$
is compact too. Recall that the homogeneous space $G/H^o$ is 
compact iff $H_\mathbb{C}$ contains a maximal connected $\mathbb{R}$-split solvable $\mathbb{R}$-subgroup of $G_\mathbb{C}$
(see, for example, \cite[Theorem 3.1]{pr}).
Since maximal connected $\mathbb{R}$-split solvable $\mathbb{R}$-subgroups of $G_\mathbb{C}$
are conjugate over $G_\mathbb{R}$ (see \cite{bt}, or Theorem 15.2.5 and Exercise 15.4.8 in \cite{spr}),
for some $g_0\in G$, $B^{g_0}\subseteq H$.
Hence, $(B^o)^{g_0}\subseteq H^o\subseteq S$.

Denote by $\pi$ the projection map $G/(B^o)^{g_0}\rightarrow G/S$.
Take 
$$
\Omega\subseteq G/S\quad\hbox{and}\quad x_0\in G/S
$$
such that $m(\partial\Omega)=0$. Set 
$$
{\Omega}^*=\pi^{-1}(\Omega)\subseteq G/(B^o)^{g_0}\quad\hbox{and}\quad {x}^*_0\in \pi^{-1}(x_0).
$$
Let $m^*$ be the $K$-invariant normalized measure on $G/(B^o)^{g_0}$. Then $m=\pi(m^*)$ is 
the $K$-invariant normalized measure on $G/S$. It is easy to check that
$\pi(\partial\Omega^*)\subseteq \partial\Omega$. Hence, 
$m^*(\partial\Omega^*)\le m^*(\pi^{-1}(\partial\Omega))=m(\partial\Omega)=0$.
Finally,
$$
N_T(\Omega,x_0)=N_T(\Omega^*,x_0^*)\mathop{\sim}_{T\rightarrow\infty}
\frac{\gamma_n}{\bar\mu(\Gamma\backslash G)} m^*(\Omega^*) T^{n(n-1)}
=\frac{\gamma_n}{\bar\mu(\Gamma\backslash G)} m(\Omega) T^{n(n-1)}.
$$

\end{proof}\medbreak


We need the following proposition that follows easily from Theorem \ref{th_ergodic}.

\begin{pro} \label{p_ergodic}
Let $f$ be the characteristic function of a relatively compact Borel subset $Z\subseteq G$
such that $\mu (\partial Z)=0$. Then for any $y\in \Gamma\backslash G$,
\begin{equation}\label{eq_char_last}
\frac{1}{\varrho(B^o_T)}\int_{B^o_T} \tilde{f}(yb^{-1})d\varrho(b)
\longrightarrow \frac{1}{\bar \mu(\Gamma\backslash G)}\int_{G} fd\mu\quad\hbox{as}\;\;T\rightarrow\infty,
\end{equation}
where $\tilde{f}(\Gamma g)=\sum_{\gamma\in \Gamma} f(\gamma g)\in C_c(\Gamma\backslash G)$.
\end{pro}

\begin{proof}[Proof.]
The argument is quite standard. One chooses functions $\phi_n,\psi_n\in C_c(G)$, $n\ge 1$,
such that $\phi_n\le f\le \psi_n$ and $\int_G (\psi_n-\phi_n)d\mu<\frac{1}{n}$.
By Theorem \ref{th_ergodic} and (\ref{eq_mubar}),
$$
\lim_{T\rightarrow\infty} \frac{1}{\varrho(B^o_T)}\int_{B^o_T} \tilde{\phi}_n(yb^{-1})d\varrho(b)
=\frac{1}{\bar \mu(\Gamma\backslash G)}\int_{\Gamma\backslash G} \tilde{\phi}_n d\bar\mu=
\frac{1}{\bar \mu(\Gamma\backslash G)}\int_{G} \phi_n d\mu,
$$
and
$$
\lim_{T\rightarrow\infty} \frac{1}{\varrho(B^o_T)}\int_{B^o_T} \tilde{\psi}_n(yb^{-1})d\varrho(b)
=\frac{1}{\bar \mu(\Gamma\backslash G)}\int_{G} \psi_n d\mu
$$
for every $n\ge 1$.
This implies (\ref{eq_char_last}).
\end{proof}\medbreak

The proof of Theorem \ref{th_asympt} should be compared with similar arguments in 
\cite{drs}, \cite{em}, \cite{ems}, \cite{emm}
where other counting problems were also reduced to asymptotics of certain integrals.

\begin{proof}[Proof of Theorem \ref{th_asympt}.]
Let $\alpha: K\rightarrow G/B^o$ be a map defined by $\alpha(k)=kB^o$.
By (\ref{eq_iwasawa}), $\alpha$ is a diffeomorphism.
The measure $m$ is given by
\begin{equation} \label{eq_mOm}
m(C)=\int_{\alpha^{-1}(C)}dk\quad\hbox{for}\;\;\hbox{Borel}\;\;\hbox{set}\;\; C\subseteq G/B^o.
\end{equation}
Since $\alpha$ is surjective, $x_0=k_0^{-1}B^o$ for some $k_0\in K$.
It follows from the Iwasawa decomposition (\ref{eq_iwasawa}) that the product map
$$
K\times B^o:(k,b)\mapsto kk_0^{-1}bk_0\in G
$$
is a diffeomorphism. For $g\in G$, define $k_g\in K$ and $b_g\in B^o$
such that 
$$
g=k_gk_0^{-1}b_gk_0.
$$
Since $G$ and $K$ are unimodular, 
it follows from (\ref{eq_iwasawa2}) that for $f\in C_c(G)$,
$$
\int_G fd\mu= \int_{K\times B^o} f(kk_0^{-1}bk_0)dk d\varrho(b).
$$

Let $\phi$ be the characteristic function of $\tilde{\Omega}\stackrel{def}{=}\alpha^{-1}(\Omega)k_0\subseteq K$,
and $\psi_{\mathcal{O}}$ the characteristic function of an open bounded symmetric
neighborhood ${\mathcal{O}}$ of $1$ in $B^o$ with boundary of measure $0$ normalized so that
$\int_{B^o} \psi_{\mathcal{O}} d\varrho=1$. 
Then $\int_{B^o} \psi_{\mathcal{O}} d\lambda=1$ too.
Note that for $g\in G$, $gx_0\in\Omega$ iff $k_g\in\tilde{\Omega}$.
Put $f_{\mathcal{O}}(g)=\phi(k_g)\psi_{\mathcal{O}}(b_g)$.
Let $\tilde{f}_{\mathcal{O}}(\Gamma g)=\sum_{\gamma\in\Gamma} f_{\mathcal{O}}(\gamma g)$.
Now Proposition \ref{p_ergodic} can be applied to $\tilde{f}_{\mathcal{O}}$:
\begin{eqnarray}
\nonumber \frac{1}{\varrho(B^o_T)}\int_{B^o_T} \tilde{f}_{\mathcal{O}} (\Gamma k_0^{-1} b^{-1} {k_0})d\varrho (b)
\mathop{\longrightarrow}_{T\rightarrow\infty} \frac{1}{\bar\mu(\Gamma\backslash G)}\int_{G} {f}_{\mathcal{O}} (gk_0)d\mu(g)\\
=\frac{1}{\bar\mu(\Gamma\backslash G)}\int_K \phi dk\cdot\int_{B^o} \psi_{\mathcal{O}} d\varrho=\frac{1}{\bar\mu(\Gamma\backslash G)}\int_{\tilde{\Omega}} dk=\frac{m(\Omega)}{\bar\mu(\Gamma\backslash G)}. \label{eq_bound0}
\end{eqnarray}
The last equality follows from (\ref{eq_mOm}).

Take $r>1$. There exists a bounded open symmetric neighborhood ${\mathcal{O}}$ 
of identity in $B^o$ (with boundary of measure $0$) such that for
any $b\in {\mathcal{O}}={\mathcal{O}}^{-1}$ and $x\in M(n,\mathbb{R})$,
\begin{equation} \label{eq_norm}
r^{-1}\|x\|\le \|b^{-1}x\|\le r\|x\|.
\end{equation}
Then for ${\mathcal{O}}$ as above,
\begin{eqnarray}
\nonumber \sum_{\gamma\in \Gamma} \int_{B^o_T} f_{\mathcal{O}} (\gamma k_0^{-1}b^{-1}{k_0})d\varrho(b)
=\sum_{\gamma\in \Gamma} \int_{(B^o_T)^{-1}} f_{\mathcal{O}} (k_\gamma k_0^{-1}b_\gamma b{k_0})d\lambda(b)\\
\nonumber =\sum_{\gamma\in \Gamma} \int_{b_\gamma(B^o_T)^{-1}} f_{\mathcal{O}} (k_\gamma k_0^{-1}b{k_0})d\lambda(b)=
\sum_{\gamma\in \Gamma} \int_{\|b^{-1}b_\gamma \|<T} \phi (k_\gamma)\psi_{\mathcal{O}}(b)d\lambda(b)\\
\nonumber =\sum_{\gamma: k_\gamma\in\tilde{\Omega}} \int_{\|b^{-1}b_\gamma \|<T} \psi_{\mathcal{O}}(b)d\lambda(b)
=\sum_{\gamma: \gamma\cdot x_0\in \Omega} \int_{\|b^{-1}b_\gamma \|<T} \psi_{\mathcal{O}}(b)d\lambda(b).
\end{eqnarray}
The integral
$$
I_\gamma\stackrel{def}{=}\int_{\|b^{-1}b_\gamma \|<T} \psi_{\mathcal{O}}(b)d\lambda(b)
$$
is not greater than $1$. By (\ref{eq_norm}), $I_\gamma=0$ for $\gamma\in\Gamma$ such that
$\|\gamma\|=\|b_\gamma\|\ge rT$. Therefore,
\begin{equation} \label{eq_bound1}
N_{rT}(\Omega, x_0)\ge\sum_{\gamma\in \Gamma} \int_{B^o_T} f_{\mathcal{O}} (\gamma k_0^{-1}b^{-1}{k_0})d\varrho(b).
\end{equation}
By (\ref{eq_norm}), $I_\gamma=1$ for $\gamma\in \Gamma$ such that $\|\gamma\|=\|b_\gamma\|< r^{-1}T$.
Thus,
\begin{equation} \label{eq_bound2}
N_{r^{-1}T}(\Omega,x_0)\le\sum_{\gamma\in \Gamma} \int_{B^o_T} f_{\mathcal{O}} (\gamma k_0^{-1}b^{-1}{k_0})d\varrho(b).
\end{equation}
It follows from (\ref{eq_BTasy}) that $\varrho(B_{r^{-1}T}^o)\sim \gamma_n r^{n-n^2} T^{n^2-n}$
as $T\rightarrow\infty$.
Then using (\ref{eq_bound1}) and (\ref{eq_bound0}), we get
\begin{eqnarray}
\nonumber \liminf_{T\rightarrow\infty} \frac{N_T(\Omega,x_0)}{T^{n^2-n}}&\ge&
\liminf_{T\rightarrow\infty} \frac{\gamma_n}{r^{n^2-n}\varrho(B^o_{r^{-1}T})}\int_{B^o_{r^{-1}T}} \tilde{f}_{\mathcal{O}} (\gamma k_0^{-1}b^{-1}{k_0})d\varrho(b)\\
\nonumber &=&\frac{\gamma_n m(\Omega)}{r^{n^2-n}\bar\mu(\Gamma\backslash G)}.
\end{eqnarray}
This inequality holds for any $r>1$. Hence,
$$
\liminf_{T\rightarrow\infty} \frac{N_T(\Omega,x_0)}{T^{n^2-n}}\ge \frac{\gamma_n m(\Omega)}{\bar\mu(\Gamma\backslash G)}.
$$
Similarly, using (\ref{eq_bound2}) and (\ref{eq_bound0}), one can prove that
$$
\limsup_{T\rightarrow\infty} \frac{N_T(\Omega,x_0)}{T^{n^2-n}}\le \frac{\gamma_n m(\Omega)}{\bar\mu(\Gamma\backslash G)}.
$$
This finishes the proof of Theorem \ref{th_asympt}.
\end{proof}\medbreak

\section{Behavior of unipotent flows} \label{sec_uni}

In this section we review some deep results on
equidistribution of unipotent flows, which are crucial for the proof of Theorem \ref{th_ergodic}.
Note that there are two different approaches available: Ratner \cite{r2} and Dani, Margulis \cite{dm93}. 
Both of these methods rely on Ratner's measure rigidity \cite{r1}. We follow the method of Dani and Margulis.
The results below were proved in \cite{dm91,dm93} for the case of one-dimensional flows
and extended to higher dimensional flows and even polynomial
trajectories in \cite{sh94,ems97,sh96}. See \cite[\S 19]{st} and \cite{kss01}
for more detailed exposition.

Appropriate adjustments are made for
the right $G$-action on $\Gamma\backslash G$ instead of left $G$-action on $G/\Gamma$.

\par{\sc Notations:}
Let $G$ be a connected semisimple Lie group
without compact factors, and $\Gamma$ a lattice in $G$.
Let $\mathfrak g$ be the Lie algebra of $G$. For positive integers $d$ and $n$, denote by
${\mathcal{P}}_{d,n}(G)$ the set of functions $q:\mathbb{R}^n\rightarrow G$ such that
for any $a,b\in\mathbb{R}^n$, the map 
$$
t\in\mathbb{R}\mapsto \hbox{Ad}(q(at+b))\in\mathfrak{g}
$$
is a polynomial of degree at most $d$ with respect to some basis of $\mathfrak{g}$.

Let $V_G=\oplus_{i=1}^{\dim \mathfrak{g}}\wedge^i\mathfrak{g}$. There is a natural
action of $G$ on $V_G$ induced from the adjoint representation. Fix a norm on $V_G$.
For a Lie subgroup $H$ of $G$ with Lie algebra $\mathfrak{h}$,  take a unit
vector $p_H\in \wedge^{\dim \mathfrak{h}}\mathfrak{h}$.

The following theorem allows us to estimate divergence of polynomial trajectories.
For its proof, see \cite[Theorems 2.1--2.2]{sh96}.
\begin{thm} \label{th_help1}
There exist closed subgroups $U_i$ ($i=1,\ldots, r$) such that each $U_i$
is the unipotent radical of a maximal parabolic subgroup, $\Gamma U_i$ is
compact in $\Gamma\backslash G$,  and for any
$d,n\in\mathbb{N}$, $\varepsilon,\delta>0$, there exists a compact set $C\subseteq \Gamma\backslash G$
such that for any $q\in {\mathcal{P}}_{d,n}(G)$ and a bounded open convex set $D\subseteq \mathbb{R}^n$,
one of the following holds:
\begin{enumerate}
\item[1.] There exist $\gamma\in\Gamma$ and $i=1,\ldots, r$ such that
$\sup_{t\in D} \|q(t)^{-1}\gamma\cdot p_{U_i}\|\le\delta.$ 
\item[2.] $\ell(t\in D: \Gamma q(t)\notin C)< \varepsilon \ell(D)$, where $\ell$
is the Lebesgue measure on $\mathbb{R}^n$.
\end{enumerate}
\end{thm}

Denote by ${\mathcal{H}}_\Gamma$ the family of all proper closed connected subgroups $H$ of $G$
such that $\Gamma\cap H$ is a lattice in $H$, and $\hbox{Ad}(H\cap \Gamma)$ is Zariski-dense in $\hbox{Ad}(H)$.
\begin{thm} \label{th_help2}
The set ${\mathcal{H}}_\Gamma$ is countable. For any $H\in {\mathcal{H}}_\Gamma$,
$\Gamma\cdot p_H$ is discrete.
\end{thm}
For proofs, see \cite[Theorem 1.1]{r1} and \cite[Theorems 2.1 and 3.4]{dm93}.

Fix a subgroup $U$ generated by $1$-parameter unipotent subgroups.
For a closed subgroup $H$ of $G$, denote
$$
X(H,U)=\{g\in G: gU\subseteq Hg\}.
$$
Define a singular set
\begin{equation} \label{eq_sing}
Y=\bigcup_{H\in {\mathcal{H}}_\Gamma} \Gamma X(H,U)\subseteq \Gamma\backslash G.
\end{equation}
It follows from Dani's generalization of Borel density theorem and
Ratner's topological rigidity that $y\in Y$ iff $yU$ is not dense in $\Gamma\backslash G$
(see \cite[Lemma 19.4]{st}).

One needs to estimate behavior of polynomial trajectories near the singular set $Y$.
The following result can be deduced from \cite[Proposition 5.4]{sh94}.
It is formulated in \cite{sh96} and \cite{kss01}.
Note that it is analogous to Theorem \ref{th_help1} with the point at infinity being
replaced by the singular set.
\begin{thm} \label{th_help3}
Let $d,n\in\mathbb{N}$, $\varepsilon>0$, $H\in {\mathcal{H}}_\Gamma$.
For any compact set $C\subseteq \Gamma X(H,U)$, there exists a compact set
$F\subseteq V_G$ such that for any neighborhood $\Phi$ of $F$ in $V_G$, there exists
a neighborhood $\Psi$ of $C$ in $\Gamma\backslash G$
such that for any $q\in {\mathcal{P}}_{d,n}(G)$ and a bounded open convex set $D\subseteq \mathbb{R}^n$,
one of the following holds:
\begin{enumerate}
\item[1.] There exists $\gamma\in\Gamma$ such that $q(D)^{-1}\gamma\cdot p_{H}\subseteq \Phi$.
\item[2.] $\ell(t\in D: \Gamma q(t)\in \Psi)< \varepsilon \ell(D)$, where $\ell$
is the Lebesgue measure on $\mathbb{R}^n$.
\end{enumerate}
\end{thm}

\section{Representations of $\hbox{\rm SL}(n,\mathbb{R})$}\label{sec_rep}

In order to be able to use the results from the previous section,
we collect here some information about representations of $\hbox{\rm SL}(n,\mathbb{R})$.

The next lemma is essentially Lemma 5.1 from \cite{sh96}.
We present its proof because more precise information about 
dependence on $\beta$ in the inequality (\ref{eq_shineq}) is needed.

\begin{lem} \label{lem_shah}
Let $V$ be a finite dimensional vector space with a norm $\|\cdot\|$,
$\mathfrak n$ be a nilpotent Lie subalgebra of $\hbox{\rm End}(V)$ with a
basis $\{b_i:i=1,\ldots, m\}$, and $N=\exp(\mathfrak{n})\subseteq \hbox{\rm GL}(V)$
be the Lie group of $\mathfrak n$.
Define a map $\Theta:\mathbb{R}^m\rightarrow N$:
$$
\Theta(t_1,\ldots, t_m)=\exp\left(\sum_{i=1}^m t_ib_i\right),\quad (t_1,\ldots, t_m)\in\mathbb{R}^m.
$$
For $\beta>0$, put $D(\beta)=\Theta\left([0,\beta]\times\cdots \times[0,\beta]\right)$. Let
$$
W=\{v\in V: \mathfrak{n}\cdot v=0\}.
$$
Denote by $\hbox{\rm pr}_W$ the orthogonal projection on W with respect to some 
scalar product on $V$.

Then there exists a constant $c_0>0$ such that for any $\beta\in (0,1)$ and $v\in V$,
\begin{equation} \label{eq_shineq}
\sup_{n\in D(\beta)} \|\hbox{\rm pr}_W(n\cdot v)\|\ge c_0\beta^d\|v\|,
\end{equation}
where $d$ is the degree of the polynomial map $\Theta$.
\end{lem}

\begin{proof}[Proof.]
Let
$$
{\mathcal{I}}=\left\{(i_1,\ldots, i_m)\in\mathbb{Z}^m: i_k\ge 0,\sum_k i_k\le d\right\}.
$$
For $t\in\mathbb{R}^m$ and $I=(i_1,\ldots, i_m)\in {\mathcal{I}}$, denote $t^I=\prod_k t_k^{i_k}$,
and $|I|=\sum_k i_k$.

One can write $\Theta(t)=\sum_{I\in {\mathcal{I}}} t^I B_I$ for some $B_I\in \hbox{End}(V)$.
Then
$$
\hbox{pr}_W(\Theta(t)v)=\sum_{I\in {\mathcal{I}}} t^I\hbox{pr}_W(B_Iv).
$$
Consider a map $T:V\rightarrow\oplus_{I\in {\mathcal{I}}} W$ defined by
$$
Tv=\sum_{I\in {\mathcal{I}}} \hbox{pr}_W(B_Iv),
$$
and a map $A_t:\oplus_{I\in {\mathcal{I}}} W\rightarrow W$ for $t\in\mathbb{R}^m$ defined by
$$
A_t\left(\bigoplus_{I\in {\mathcal{I}}} w_I\right)=\sum_{I\in {\mathcal{I}}} t^Iw_I.
$$
For ${I\in {\mathcal{I}}}$, take fixed $s_I\in\mathbb{R}^m$ such that $0<s_{I,k}<1$ and $s_{I_1,k}\ne s_{I_2,k}$
for $I_1\ne I_2$, and put $t_I=\beta s_I$. Let 
$$
A=\bigoplus_{I\in {\mathcal{I}}} A_{t_I}:\bigoplus_{I\in {\mathcal{I}}} W\rightarrow\bigoplus_{I\in {\mathcal{I}}} W.
$$
The map $A$ has a matrix form $\left(t_I^J\right)_{I,J\in {\mathcal{I}}}$. This matrix is a Kronecker
product of Vandermonde matrices which implies that $A$ is invertible. Using elementary row and column
operations, one can write
\begin{equation} \label{eq_norm1}
\left(t_I^J\right)_{I,J\in {\mathcal{I}}}=B\cdot\hbox{diag}\left(\beta^{|I|}:I\in {\mathcal{I}}\right)\cdot C
\end{equation}
for some $B,C\in\hbox{GL}\left(\oplus_{I\in {\mathcal{I}}} W\right)$, which are independent of $\beta$.
It is convenient to use a norm on $\oplus_{I\in {\mathcal{I}}} W$ defined by
$$
\left\|\bigoplus_{I\in {\mathcal{I}}} w_I\right\|=\max_{I\in {\mathcal{I}}} \|w_I\|,\quad w_I\in W.
$$
Then by (\ref{eq_norm1}), for $w\in\oplus_{I\in {\mathcal{I}}} W$,
\begin{equation} \label{eq_norm2}
\|Aw\|\ge \|B^{-1}\|^{-1}\cdot\beta^d\cdot\|C^{-1}\|^{-1}\cdot\|w\|.
\end{equation}

It follows from Lie-Kolchin theorem that $T$ is injective (see \cite[Lemma 5.1]{sh96}). Therefore, there
exists a constant $c_1>0$ such that $\|Tv\|\ge c_1\|v\|$ for $v\in V$.
Then using (\ref{eq_norm2}), we get
$$
\sup_{n\in D(\beta)} \|\hbox{\rm pr}_W(n\cdot v)\|=\sup_{t:0<t_i<\beta} \|A_tTv\|\ge \|ATv\|\ge c_0\beta^d\|v\|,
$$
where $c_0=\|B^{-1}\|^{-1}\cdot\|C^{-1}\|^{-1}\cdot c_1$.
\end{proof}\medbreak

We will need the following elementary observation:

\begin{lem} \label{lem_polynom}
Let $F:\mathbb{R}^m\rightarrow\mathbb{R}^m$ be a $C^1$ bijection such that $F(0)=0$,
and $F^{-1}$ is $C^1$ too.
Fix a norm on $\mathbb{R}^m$ and denote by $D(r)$ a ball of radius $r$ centered at the origin. 
Then there are $c_1,c_2>0$ such that
$$
D(c_1r)\subseteq F(D(r))\subseteq D(c_2r)
$$
for every $r\in (0,1)$. 
\end{lem}



Let $\mathfrak g$ be the Lie algebra of $G=\hbox{\rm SL}(n,\mathbb{R})$,
and $\mathfrak{g}_\mathbb{C}=\mathfrak{g}\otimes \mathbb{C}$.
Recall the root space decomposition of $\mathfrak{g}_\mathbb{C}$:
$$
\mathfrak{g}_\mathbb{C}=\mathfrak{g}_0\oplus \sum_{i\ne j} \mathfrak{g}_{ij},
$$
where $\mathfrak{g}_0$ is the subalgebra of diagonal matrices of $\mathfrak{g}_\mathbb{C}$,
and $\mathfrak{g}_{ij}=\mathbb{C}E_{ij}$ ($E_{ij}$ is the matrix with $1$
in position $(i,j)$ and $0$'s elsewhere).

Introduce {\it fundamental weights} of $\mathfrak{g}_\mathbb{C}$:
\begin{equation} \label{eq_fw}
\omega_i (s)=s_1+\cdots +s_i,\quad 1\le i\le n-1,
\end{equation}
where $s\in \mathbb{C}^n$ and $\sum_i s_i=0$.
{\it Dominant weights} are defined as linear combinations with non-negative integer coefficients of
the fundamental weights $\omega_i$, $1\le i\le n-1$.

A {\it highest weight} of a representation of $\mathfrak{g}_\mathbb{C}$ is a weight
that is maximal with respect to the ordering on the dual space of $\mathfrak{g}_0$.
Recall that irreducible representations of $\mathfrak{g}_\mathbb{C}$ are classified by their highest weights
(see, for example, \cite[Ch.~5]{gw}). 
The highest weights are precisely the dominant weights
defined above.

\begin{lem} \label{lem_diverge0}
Let $\pi$ be a representation of $G$ on a finite dimensional complex vector space $V$.
Let $x\in V-\{0\}$ be such that $\pi(N)x=x$.
Then $x$ is a sum of weight vectors with dominant weights.
Moreover, if $V$ does not contain non-zero $G$-fixed vectors, every weight
in this sum is not zero.
\end{lem}

\begin{proof}[Proof.]
Consider a representation $\tilde \pi$ of $\mathfrak{g}_\mathbb{C}$ on $V$ induced by the representation $\pi$. 
Since this representation is completely reducible, it is enough to
consider the case when it is irreducible.

We claim that in this case, $x$ is a vector with the highest weight. Write $x=\sum_k x_k$ where each $x_k\in V$
is a weight vector with a weight $\lambda_k$. We may assume that $\lambda_k\ne \lambda_l$
for $k\ne l$. Since $\pi(N)x=x$, $\tilde\pi(E_{ij})x=0$ for $i<j$. Thus,
$\sum_k \tilde\pi(E_{ij})x_k=0$. Since $\tilde\pi(E_{ij})x_k$ is either $0$ or
a weight vector with the weight $\lambda_k+\alpha_{ij}$, the non-zero terms in the sum
are linearly independent. Hence, $\tilde\pi(E_{ij})x_k=0$ for $i<j$.
Suppose that $\lambda_k$ is not the highest weight. Note that $\tilde\pi(\mathfrak{g}_0)x_k=\mathbb{C}x_k$,
and $\tilde\pi(E_{ji})x_k$ has weight $\lambda_k-\alpha_{ij}<\lambda_k$ for $i<j$.
By Poincar\'e-Birkhoff-Witt theorem, the universal enveloping algebra
$\mathcal{U}(\mathfrak g_\mathbb{C})=\mathcal{U}(\mathfrak b^-)\oplus\mathcal{U}(\mathfrak g_\mathbb{C})\mathfrak n$,
where $\mathfrak b^-$ is the space of lower triangular matrices, and $\mathfrak n$ is the Lie algebra of $N$.
Therefore, the space $\tilde\pi(\mathcal{U}(\mathfrak{g}_\mathbb{C}))x_k=\tilde\pi(\mathcal{U}(\mathfrak{b}^-))x_k$ does not contain a vector with the highest weight.
This contradicts irreducibility of $\tilde\pi$. Thus, each $x_k$ is of the highest weight,
and $x$ is a highest weight vector.
Since every highest weight is dominant, the lemma is proved.
\end{proof}\medbreak

For a fixed $g_0\in G$, define $q_s(t)=g_0n(t)^{-1}a(s)^{-1}$. We are going to study $q_s$ using techniques from Section
\ref{sec_uni}. The next lemma guarantees that the first possibility in
Theorems \ref{th_help1} and \ref{th_help3} does not occur.

For $\beta>0$, define
\begin{equation} \label{eq_drho}
D(\beta)=\left\{n(t)\in N: \sum_{i<j} t_{ij}^2<\beta^2\right\}.
\end{equation}

\begin{lem} \label{l_diverge}
Let $\pi$ be a nontrivial representation of $G$ on a real vector space $V$ equipped with a norm $\|\cdot\|$.
Let $V_0=\{v\in V: \pi(G)v=v\}$, and $V_1$ be a $G$-invariant complement for $V_0$.
(Note that $V_1$ exists because $\pi$ is completely reducible.)
Denote by $\Pi$ the projection on $V_1$.

Then for any relatively compact set $K\subseteq V$ and
$r>0$, there exist $\alpha \in (0,1)$ and $C_0>0$ such that for any $s$ with $a(s)\in A_T^{C_0}$ and
$x\in V$ with $\|\Pi(x)\|>r$,
\begin{equation}\label{eq_div}
\pi\left(q_s\left(D(e^{-\alpha s_1})\right)\right)^{-1}\cdot x\nsubseteq K.
\end{equation}
\end{lem}

\begin{proof}[Proof.]
It is convenient to extend $\pi$ to $V_\mathbb{C}=V\otimes \mathbb{C}$.
$(V_0)_\mathbb{C}$ is the space of $G$-fixed vectors in $V_\mathbb{C}$.
Thus, we may assume $V$ to be complex. Also dealing with projections on $V_1$, we may assume
that $V$ has no $G$-fixed vectors.

Since $\{g_0^{-1}\cdot x:\|x\|>r\}\subseteq \{x:\|x\|>r_1\}$ for some $r_1>0$,
we may assume that $g_0=1$.

For some $R>0$, $K\subseteq \{x\in V: \|x\|<R\}$. If (\ref{eq_div})
fails for some $s\in\mathbb{R}^n$ and $x\in V$, then 
\begin{equation} \label{eq_g0}
\sup_{n\in D(e^{-\alpha s_1})} \|\pi(a(s)n)x\|<R.
\end{equation}

Let $W=\{x\in V:\pi(N)x=x\}$.
Clearly, the statement of the lemma is independent of the norm used. It is convenient
to use the max-norm with respect to a basis $\left\{ v_i\right\}$ of $V$ consisting of weight
vectors, i.e.
$$
\left\|\sum_i u_i v_i\right\|=\max_i |u_i|,\quad u_i\in \mathbb{C},
$$
and each $v_i$ is an eigenvector of $A$.
Moreover we can choose the basis $\{v_i\}$ so that it contains a basis of $W$.
Let $\hbox{pr}_W$ be the projection onto $W$ with respect to this basis.
Then $\hbox{pr}_W$ commutes with $a(s)$.
Note that there exists $C'>0$ such that for any $v\in V$,
\begin{equation} \label{eq_g1}
\|v\|\ge C'\|\hbox{pr}_W(v)\|.
\end{equation}

Let $\mathcal{K}\subseteq \mathbb{N}$ be such that $k\in \mathcal{K}$
iff $v_k$ has a non-zero dominant weight. Denote this weight by $\lambda_k$.
By Lemma \ref{lem_diverge0}, $W\subseteq \left<v_k:k\in \mathcal{K}\right>$.
In particular, for $n\in N$, 
\begin{equation} \label{eq_g2}
\hbox{pr}_W(\pi(n)x)=\sum_{k\in {\mathcal{K}}} u_k(n)v_k\quad\hbox{for}\;\hbox{some}\; u_k(n)\in \mathbb{C}.
\end{equation}
Therefore, 
\begin{equation} \label{eq_g3}
\|\hbox{pr}_W(\pi(n)x)\|=\max_{k\in {\mathcal{K}}} |u_k(n)|.
\end{equation}
Using the fact that $\hbox{pr}_W$ and $\pi(a(s))$ commute, 
(\ref{eq_g1}), (\ref{eq_g2}), and (\ref{eq_g3}), we have that for any $n\in N$,
\begin{eqnarray}
\|\pi(a(s)n)x\|\ge C'\|\hbox{pr}_W(\pi(a(s))\pi(n)x)\|= 
C'\left\|\pi(a(s))\left(\sum_{k\in {\mathcal{K}}} u_k(n)v_k\right)\right\| \nonumber\\
=C'\max_{k\in {\mathcal{K}}} \left(|u_k(n)|e^{\lambda_k(s)}\right)\ge C'\exp\left(\min_{k\in {\mathcal{K}}} \lambda_k(s)\right)\|\hbox{pr}_W(\pi(n)x)\|. \label{eq_g4}
\end{eqnarray}

Let $\mathfrak n$ be the Lie algebra of $N$. 
Denote by $\tilde\pi$ the representation of $\mathfrak g$ induced by $\pi$.
Since $\mathfrak g$ is simple, $\tilde\pi$ is faithful.
Thus, $\tilde\pi$ defines an isomorphism between $\mathfrak n$ and $\tilde\pi(\mathfrak n)$.
Since the exponential map $\mathfrak n\to N$ is a polynomial isomorphism,
the coordinates on $N$ used in Lemma \ref{lem_shah}
and the coordinates $\{t_{ij}\}$ are connected by a polynomial isomorphism too.
By Lemma \ref{lem_polynom}, (\ref{eq_shineq}) holds for the set $D(\beta)$ defined in (\ref{eq_drho}).
Therefore,
\begin{equation}\label{eq_g5}
\sup_{n\in D(e^{-\alpha s_1})}\|\hbox{pr}_W(\pi(n)x)\|\ge c_0 (e^{-\alpha s_1})^d\|x\|\ge c_0re^{-\alpha ds_1}
\end{equation}
for some positive integer $d$.
It follows from (\ref{eq_g4}) and (\ref{eq_g5}) that if (\ref{eq_g0}) holds, then
$$
\exp\left(\min_{k\in {\mathcal{K}}} \lambda_k (s)-\alpha d s_1\right)\le \frac{R}{c_0C'r}.
$$
Take $\alpha<d^{-1}$. Since each $\lambda_k$ is a non-zero dominant weight, it follows from
(\ref{eq_fw}) that that $\lambda_k(s)-\alpha d s_1\rightarrow\infty$
as $C\rightarrow\infty$ for $s$ such that $a(s)\in A^C$.
Hence, there exists $C_0>0$ such that (\ref{eq_g0}) does not hold
for $s$ with $a(s)\in A_T^{C_0}$.
Since (\ref{eq_g0}) fails, (\ref{eq_div}) holds.
\end{proof}\medbreak

\begin{lem} \label{l_diverge00}
Use notations from Lemma \ref{l_diverge}.
Let $H$ be a subgroup of $G$, and $x\in V$ such that $\pi(H)x$ is discrete in $V$.
Then $\Pi(\pi(H)x)$ is discrete in $V_1$.
\end{lem}

\begin{proof}[Proof.]
Denote by $x_0\in V_0$ and $x_1\in V_1$ the components of $x$ with respect to
the decomposition $V=V_0\oplus V_1$. Then $\Pi(\pi(H)x)=\pi(H)x_1$.
Suppose that for some $\{h_n\}\subseteq H$, $\pi(h_n)x_1\rightarrow y$
for some $y\in V_1$. Then $\pi(h_n)x$ converges to $x_0+y$. It follows that
$\pi(h_n)x$ is constant for large $n$. Therefore, $\pi(h_n)x_1=\pi(h_n)x-x_0$
is constant for large $n$ too.
\end{proof}\medbreak

\section{Proof of Theorem \ref{th_ergodic}}\label{sec_th2}

Let ${\mathcal{Z}}=(\Gamma\backslash G)\cup\{\infty\}$ be the $1$-point compactification
of $\Gamma\backslash G$.
For $T>0$, define a normalized measure on ${\mathcal{Z}}$ by
$$
\nu_T(\tilde{f})=\frac{1}{\varrho(B^o_T)}\int_{B^o_T} \tilde{f}(yb^{-1})d\varrho(b),\quad \tilde{f}\in C_c(\Gamma\backslash G).
$$
To prove Theorem \ref{th_asympt}, 
we need to show that $\nu_T\rightarrow \nu$ in weak-$*$ topology as $T\rightarrow\infty$.
Since the space of normalized measures on ${\mathcal{Z}}$ is compact,
it is enough to prove that if $\nu_{T_i}\rightarrow \eta$ as $T_i\rightarrow\infty$ for some normalized
measure $\eta$ on ${\mathcal{Z}}$, then $\eta$ is $G$-invariant, and $\eta (\{\infty\})=0$.

It follows from Lemma \ref{lem_BTC2} that for any $C\in\mathbb{R}$,
\begin{equation} \label{eq_eta}
\eta(\tilde{f})=\lim_{T_i\rightarrow\infty}
\frac{1}{\varrho(B^o_{T_i})}\int_{B^C_{T_i}} \tilde{f}(yb^{-1})d\varrho(b),\quad \tilde{f}\in C_c(\Gamma\backslash G).
\end{equation}

Let
\begin{equation}\label{eq_U}
U=\{n(t)\in N: t_{ij}=0\;\hbox{for}\; i<j<n\}.
\end{equation}

\begin{lem} \label{lem_eta}
The measure $\eta$ is $U$-invariant.
\end{lem}

\begin{proof}[Proof.]
For $\tilde{f}\in C_c(\Gamma\backslash G)$, and $g_0\in G$, define 
$\tilde{f}^{g_0}(\Gamma g)=\tilde{f}(\Gamma gg_0)\in C_c(\Gamma\backslash G)$.

Let $\tilde{f}\in C_c(\Gamma\backslash G)$. Take $M>0$ such that $|\tilde{f}|<M$.

For $T>0$ and $s\in\mathbb{R}^{n-1}$, define a set 
\begin{equation} \label{eq_Dts}
D_{s,T}=\{n\in N: \|a(s)n\|<T\}.
\end{equation}
Denote by $\chi_{s,T}(n)$ the characteristic function of the set $D_{s,T}$.
Then we can rewrite (\ref{eq_eta}) as
\begin{equation}\label{eq_eta1}
\eta(\tilde{f})=\lim_{T_i\rightarrow\infty}\frac{1}{\varrho(B^o_{T_i})}
\int_{A_{T_i}^C}\int_{N} \tilde{f}(yn^{-1}a(-s))\chi_{s,T_i}(n)e^{2\delta(s)}dnds.
\end{equation}
Let 
\begin{equation}\label{eq_apr}
A_T^{'C}=\{a(s)\in A_T^C: T^2-N(s)> T\},
\end{equation}
where $N(s)$ is defined in (\ref{eq_Ns}), and $B_T^{'C}=A_T^{'C}N\cap B_T^o$.
We claim that the equality (\ref{eq_eta1})
holds when $A_{T_i}^C$ is replaced by $A_{T_i}^{'C}$. By Lemma \ref{lem_BTC},
\begin{eqnarray*}
&&\int_{A_{T_i}^C-A_{T_i}^{'C}}\int_{N} \tilde{f}(yn^{-1}a(-s))\chi_{s,T_i}(n)e^{2\delta(s)}dnds\\
&\ll& \int_{A_{T_i}^C-A_{T_i}^{'C}}\left(T_i^2-N(s)\right)^{\frac{n(n-1)}{4}}\exp\left(\sum_k (n-k)s_k\right)ds\\
&=& O\left(T_i^{\frac{3n(n-1)}{4}}\right)\quad\hbox{as}\quad T_i\rightarrow\infty.
\end{eqnarray*}
Now the claim follows from (\ref{eq_BTasy}).

Take $u\in U$.
Let $u(s)=\hbox{Ad}_{a(-s)}(u)$. Then
\begin{eqnarray}
\nonumber |\eta(\tilde{f}^{u})-\eta(\tilde{f})|\le
\limsup_{T_i\rightarrow\infty}\frac{1}{\varrho(B_{T_i}^o)}\int_{A_{T_i}^{'C}}\int_{N} |\tilde{f}(yn^{-1}u(s)a(-s))\chi_{s,T_i}(n)\\
- \tilde{f}(yn^{-1}a(-s))\chi_{s,T_i}(n)|dn e^{2\delta(s)}ds. \label{eq_mess1}
\end{eqnarray}
We estimate the last integral:
\begin{eqnarray}
\nonumber \int_{N} \left|\tilde{f}(yn^{-1}u(s)a(-s))\chi_{s,T_i}(n)-\tilde{f}(yn^{-1}a(-s))\chi_{s,T_i}(n)\right|dn\\
\nonumber =\int_{N} \left|\tilde{f}(yn^{-1}a(-s))\right|\cdot \left|\chi_{s,T_i}(u(s)n)-\chi_{s,T_i}(n)\right|dn\\
\le M \int_{N} \left|\chi_{s,T_i}(u(s)n)-\chi_{s,T_i}(n)\right|dn. \label{eq_mess2}
\end{eqnarray}

Recall that $\alpha_{i,n}(-s)=-s_i+s_n$. Therefore, by (\ref{eq_adj}),
$u(s)_{in}= e^{-s_i+s_n} u_{in}$ for $i=1,\ldots, n-1$. It follows from the triangle
inequality that 
$$
\|a(s)u(s)n\|\le \|a(s)n\|+\sqrt{\sum_{i=1}^{n-1} e^{2s_i}u(s)_{in}^2}\le \|a(s)n\|+ e^{s_n}\|u\|,
$$
and similarly,
$$
\|a(s)n\|\le \|a(s)u(s)n\|+ e^{s_n}\|u^{-1}\|=\|a(s)u(s)n\|+ e^{s_n}\|u\|.
$$
Hence,
$$
\chi_{s,T_i-e^{s_n}\|u\|}(n)\le \chi_{s,T_i}(u(s)n)\le \chi_{s,T_i+e^{s_n}\|u\|}(n)
$$
for $n\in N$. Therefore,
\begin{equation}\label{eq_mess3}
\int_{N} \left|\chi_{s,T_i}(u(s)n)-\chi_{s,T_i}(n)\right|dn\le\int_{N}\left(\chi_{s,T_i+e^{s_n}\|u\|}-
\chi_{s,T_i-e^{s_n}\|u\|}\right)dn.
\end{equation}

Let $\varepsilon>0$. We claim that there exists $C_0>0$ such that for $C>C_0$ and 
$a(s)\in A^{'C}_{T_i}$,
\begin{equation} \label{eq_mess4}
\int_{N}\left(\chi_{s,T_i+e^{s_n}\|u\|}-\chi_{s,T_i-e^{s_n}\|u\|}\right)dn
\le \varepsilon \int_{N}\chi_{s,T_i}dn.
\end{equation}
Similarly to Lemma \ref{lem_BTC},
$$
\int_{N}\chi_{s,T_i}dn=c_n\left(T_i^2-N(s)\right)^{\frac{n(n-1)}{4}}\exp\left(\sum_k (n-k)s_k\right).
$$
Also $e^{s_n}\rightarrow 0$ for $a(s)\in A^{'C}_{T_i}$ as $C\rightarrow\infty$.
Therefore, the equation (\ref{eq_mess4}) will follow from the following.

{\it Claim}. There exists $d_0=d_0(\varepsilon)>0$ such that for any $d\in (0,d_0)$ and
$a(s)\in A^{'C}_{T_i}$,
\begin{equation}\label{eq_LLL}
\Big((T_i+d)^2-N(s)\Big)^{\frac{n(n-1)}{4}}-\Big((T_i-d)^2-N(s)\Big)^{\frac{n(n-1)}{4}}
<\varepsilon \Big((T_i-d)^2-N(s)\Big)^{\frac{n(n-1)}{4}}.
\end{equation}

Note that by (\ref{eq_apr}), $(T_i-d)^2-N(s)>0$ for $a(s)\in A^{'C}_{T_i}$ and $d<1/2$.
The inequality (\ref{eq_LLL}) is equivalent to
$$
(T_i+d)^2-N(s)<(1+\bar\varepsilon)\Big((T_i-d)^2-N(s)\Big),
$$
where $\bar\varepsilon=(1+\varepsilon)^\frac{4}{n(n-1)}-1$.
By (\ref{eq_apr}), the last inequality follows from
$$
(T_i+d)^2 +\bar\varepsilon (T_i^2-T_i)<(1+\bar\varepsilon)(T_i-d)^2,
$$
or equivalently,
$$
T_i(4d-\bar\varepsilon+2d\bar\varepsilon)<\bar\varepsilon d^2.
$$
If one takes $d<d_0=\bar\varepsilon/(2\bar\varepsilon+4)$, the left hand side is negative.
This proves the claim (\ref{eq_LLL}).

Thus, (\ref{eq_mess4}) holds.
Then by (\ref{eq_mess1}), (\ref{eq_mess2}), (\ref{eq_mess3}), and (\ref{eq_mess4}),
\begin{eqnarray*}
|\eta(\tilde{f}^{u})-\eta(\tilde{f})|&\le&
(M\varepsilon)\limsup_{T_i\rightarrow\infty}\frac{1}{\varrho(B_{T_i}^o)}\int_{A_{T_i}^{'C}}\int_{N} \chi_{s,T_i}(n)e^{2\delta(s)} dnds\\
&=&(M\varepsilon)\limsup_{T_i\rightarrow\infty}\frac{\varrho(B_{T_i}^{'C})}{\varrho(B_{T_i}^o)}\le M\varepsilon.
\end{eqnarray*}
This shows that $\eta(\tilde{f}^{u})=\eta(\tilde{f})$.
\end{proof}\medbreak

\begin{lem} \label{lem_tildeA}
Let $\alpha \in (0,1)$. Let 
\begin{equation} \label{eq_tildeA}
\tilde{A}^C_T=\left\{a(s)\in A_T^C|\; (T^2-N(s))^{1/2}>\exp\left(\mathop{\max}_{1\le j\le n-1}\{s_j\}-\alpha s_1\right)\right\},
\end{equation}
where $N(s)$ is as in (\ref{eq_Ns}), and $\tilde{B}^C_T=\tilde{A}^C_TN\cap B_T^o$. Then for $C>0$,
$$
\eta(\tilde{f})=\lim_{T_i\rightarrow\infty}
\frac{1}{\varrho(B^o_{T_i})}\int_{\tilde{B}^C_{T_i}} \tilde{f}(yb^{-1})d\varrho(b),\quad \tilde{f}\in C_c(\Gamma\backslash G).
$$
\end{lem}

\begin{proof}[Proof.]
By (\ref{eq_eta}), it is enough to show that
\begin{equation} \label{eq_tBTC}
\frac{\varrho\left(B^C_{T_i}-\tilde{B}^C_{T_i}\right)}{\varrho\left(B^o_{T_i}\right)}\rightarrow 0\quad\hbox{as}\quad T_i\rightarrow \infty.
\end{equation}
As in Lemma \ref{lem_BTC},
$$
\varrho\left(B_{T_i}^C-\tilde{B}^C_{T_i}\right)=c_n\mathop{\int}_{A_{T_i}^C-\tilde{A}^C_{T_i}} \Big(T_i^2-N(s)\Big)^{\frac{n(n-1)}{4}}\hbox{\rm exp}\left(\sum_{k} (n-k)s_k\right)ds.
$$
Therefore,
\begin{eqnarray*}
\varrho\left(B_{T_i}^C-\tilde{B}^C_{T_i}\right)\le c_n \int_{A_{T_i}^C}\exp\Big(\frac{n(n-1)}{2}\mathop{\max}_{1\le j\le n-1}\{s_j\}\\
-\frac{\alpha n(n-1)}{2} s_1+\sum_{k} (n-k)s_k\Big)ds
\le c_n \mathop{\sum}_{1\le j\le n-1} \int_{A_{T_i}^C}\exp\Big(\frac{n(n-1)}{2}s_j\\
-\frac{\alpha n(n-1)}{2} s_1+\sum_{k} (n-k)s_k\Big)ds
\end{eqnarray*}
Then as in the proof of Lemma \ref{lem_BTC2}, for $j\ne 1$
\begin{eqnarray*}
\int_{A_{T_i}^C}\exp\left(\frac{n(n-1)}{2}s_j-\frac{\alpha n(n-1)}{2} s_1+\sum_{k} (n-k)s_k\right)ds\\
\le \int_{-\infty}^{\log T_i} \exp\left(\left(-\frac{\alpha n(n-1)}{2}+n-1\right)s_1\right)ds_1\\
\cdot\int_{-\infty}^{\log T_i} \exp\left(\left(\frac{n(n-1)}{2}+n-j\right)s_j\right)ds_j\cdot
\prod_{{k<n},{k\ne 1,j}} \int_{-\infty}^{\log T_i} e^{(n-k)s_k}ds_k\\
\ll T_i^{n(n-1)-\alpha n(n-1)/2}
\end{eqnarray*}
as $T_i\rightarrow\infty$.
For $j=1$, the same estimate can be obtained by a similar calculation.
Now (\ref{eq_tBTC}) follows from (\ref{eq_BTasy}).
\end{proof}\medbreak

Let $y=\Gamma g_0$ for $g_0\in G$. Define $q_s(t)=g_0n(t)^{-1}a(s)^{-1}$.
We apply the results of Section \ref{sec_uni} to the map $q_s$.

\begin{lem} \label{lem_inf1}
$\eta(\{\infty\})=0$.
\end{lem}

\begin{proof}[Proof.]
Let $\varepsilon,\delta >0$. Apply Theorem \ref{th_help1} to the map $q_s(t)$.
By Theorem \ref{th_help2}, the set $\Gamma\cdot p_{U_i}\subseteq V_G$ is discrete.
Write $V_G=V_0\oplus V_1$, where $V_0$ is the space of vectors fixed by $G$,
and $V_1$ is its $G$-invariant complement. Denote by $\Pi$ the projection on $V_1$.
By Lemma \ref{l_diverge00}, $\Pi(\Gamma\cdot p_{U_i})$ is discrete.
Also $0\notin \Pi(\Gamma\cdot p_{U_i})$. Otherwise $p_{U_i}$ is fixed by $G$, and
it would follow that $U_i$ is normal in $G$ which is a contradiction.
Therefore, there exists $r>0$ such that 
$$
\|\Pi(x)\|>r\quad\textrm{for}\quad x\in\bigcup_i \Gamma\cdot p_{U_i}.
$$
Now we can apply Lemma \ref{l_diverge}. Let 
$$
K=\{x\in V_G: \|x\|\le \delta\}.
$$
By Lemma \ref{l_diverge},
there exist $\alpha \in (0,1)$ and $C_0>0$ such that the first case of Theorem \ref{th_help1}
fails for $q_s$ when $D$ is a bounded open convex set which contains
$D(e^{-\alpha s_1})$ (it is defined in (\ref{eq_drho})),
and $s$ is such that $a(s)\in A^{C_0}_{T_i}$. Therefore, for some compact set
$C\subseteq \Gamma\backslash G$,
\begin{equation} \label{eq_mC}
\omega \left(\left\{n(t)\in D: \Gamma q_s(t)\notin C\right\}\right)< \varepsilon
\omega (D)
\end{equation}
when $D\supseteq D(e^{-\alpha s_1})$ and $a(s)\in A^{C_0}_{T_i}$,
where $\omega=dt$ denotes the Lebesgue measure on $N$.

Let $D_{s,T_i}$ be as in (\ref{eq_Dts}). By Lemma \ref{lem_tildeA},
\begin{equation} \label{eq_tA}
\eta(\tilde{f})=\lim_{T_i\rightarrow\infty}
\frac{1}{\varrho(B^o_{T_i})}\int_{\tilde{A}^C_{T_i}}\int_{D_{s,T_i}}\tilde{f}(\Gamma q_s(t))e^{2\delta(s)}dtds,\quad \tilde{f}\in C_c(\Gamma\backslash G).
\end{equation}
Note that
$$
D_{s,T_i}=\left\{n(t)\in N: \sum_{i<j} e^{2s_i}t_{ij}^2<T_i^2-N(s)\right\},
$$
where $N(s)$ is defined in (\ref{eq_Ns}) (cf. (\ref{eq_BT})).
It follows that $D_{s,T_i}$ contains $D(\beta)$ for 
$$
\beta<(T_i^2-N(s))^{1/2}\exp\left(-\mathop{\max}_{1\le i\le n-1} \{s_i\}\right).
$$
When $a(s)\in\tilde{A}^{C_0}_{T_i}$, the right hand side is greater then $e^{-\alpha s_1}$
(see (\ref{eq_tildeA})). Therefore, $D_{s,T_i}\supseteq D(e^{-\alpha s_1})$ when
$a(s)\in\tilde{A}^{C_0}_{T_i}$. By (\ref{eq_mC}),
\begin{equation} \label{eq_mC1}
\omega \left(\left\{n(t)\in D_{s,T_i}: \Gamma q_s(t)\notin C\right\}\right)< \varepsilon
\omega (D_{s,T_i})\quad \hbox{for}\;\; a(s)\in\tilde{A}^{C_0}_{T_i}.
\end{equation}

Let $\chi_C$ be the characteristic function of the set $C$.
Take $\tilde{f}\in C_c(\Gamma\backslash G)$
such that $\chi_C\le \tilde{f}\le 1$. Then using (\ref{eq_tA}) and (\ref{eq_mC1}), we get
\begin{eqnarray*}
\eta(\hbox{supp}(\tilde{f}))&\ge& \lim_{T_i\rightarrow\infty}
\frac{1}{\varrho(B^o_{T_i})}\int_{\tilde{A}^{C_0}_{T_i}}\int_{D_{s,T_i}}\chi_C(\Gamma q_s(t))e^{2\delta(s)}dtds\\
&\ge&\lim_{T_i\rightarrow\infty}
\frac{1}{\varrho(B^o_{T_i})}\int_{\tilde{A}^{C_0}_{T_i}}(1-\varepsilon)\omega(D_{s,T_i})e^{2\delta(s)}ds\\
&=&(1-\varepsilon)\lim_{T_i\rightarrow\infty} \frac{\varrho(\tilde{B}^{C_0}_{T_i})}{\varrho(B^o_{T_i})}=1-\varepsilon.
\end{eqnarray*}
Hence, $\eta(\{\infty\})\le\eta(\hbox{supp}(\tilde{f})^c)\le\varepsilon$ for every $\varepsilon>0$.
\end{proof}\medbreak

Recall that the singular set $Y$ of $U$ was defined in (\ref{eq_sing}).

\begin{lem} \label{lem_Y1}
$\eta(Y)=0$.
\end{lem}

\begin{proof}[Proof.]
By (\ref{eq_sing}) and Theorem \ref{th_help2}, it is enough to show that
$\eta(\Gamma X(H,U))=0$ for any $H\in {\mathcal{H}}_\Gamma$. Moreover, since $\Gamma X(H,U)$
is $\sigma$-compact, we just need to show that $\eta(C)=0$ for any compact set
$C\subseteq \Gamma X(H,U)$.

Take $\varepsilon>0$. Apply Theorem \ref{th_help3} to the map $q_s(t)$.
Let $F\subseteq V_G$ be as in Theorem \ref{th_help3}.
Fix a relatively compact neighborhood $\Phi$ of $F$ in $V_G$.
Take $\Psi\supseteq C$ as in Theorem \ref{th_help3}.
By Theorem \ref{th_help2}, the set $\Gamma\cdot p_H$ is discrete.
Let $\Pi$ be as in the proof of Lemma \ref{lem_inf1}.
By Lemma \ref{l_diverge00}, $\Pi(\Gamma\cdot p_H)$ is discrete.
If $0\in\Pi(\Gamma\cdot p_H)$, the vector $p_H$ is fixed by $G$, and $H$
is normal in $G$, which is impossible.
Therefore, for some $r>0$, $\|\Pi(x)\|>r$ for every $x\in \Gamma\cdot p_H$.
Applying Lemma \ref{l_diverge} with $K=\Phi$, one gets that
there exist $\alpha\in (0,1)$ and $C_0>0$ such that the first case of Theorem \ref{th_help1}
fails for $q_s$ when $D$ is a bounded open convex set containing
$D(e^{-\alpha s_1})$, and $s$ is such that $a(s)\in A^{C_0}_{T_i}$. Therefore, 
the second case should hold:
\begin{equation} \label{eq_mC3}
\omega \left(\left\{n(t)\in D: \Gamma q_s(t)\in \Psi\right\}\right)< \varepsilon \omega (D)
\end{equation}
when $D\supseteq D(e^{-\alpha s_1})$ and $a(s)\in A^{C_0}_{T_i}$.

Let $D_{s,T_i}$ be as in (\ref{eq_Dts}). It was shown in the proof of Lemma \ref{lem_inf1}
that $D_{s,T_i}\supseteq D(e^{-\alpha s_1})$ when $a(s)\in\tilde{A}^{C_0}_{T_i}$.
It follows from (\ref{eq_mC3}) that
\begin{equation} \label{eq_mC4}
\omega \left(\left\{n(t)\in D_{s,T_i}: \Gamma q_s(t)\in \Psi\right\}\right)< \varepsilon
\omega (D_{s,T_i})\quad \hbox{for}\;\; a(s)\in\tilde{A}^{C_0}_{T_i}.
\end{equation}

Take a function $\tilde{f}\in C_c(\Gamma\backslash G)$ such that $\tilde{f}=1$ on $C$,
$\hbox{supp}(\tilde{f})\subseteq \Psi$, and $0\le \tilde{f}\le 1$.
Let $\chi_\Psi$ be the characteristic function of $\Psi$.
Then using (\ref{eq_tA}) and (\ref{eq_mC4}), we get
\begin{eqnarray*}
\eta(C)\le \lim_{T_i\rightarrow\infty}
\frac{1}{\varrho(B^o_{T_i})}\int_{\tilde{A}^{C_0}_{T_i}}\int_{D_{s,T_i}}\chi_\Psi(\Gamma q_s(t))e^{2\delta(s)}dtds\\
\le \lim_{T_i\rightarrow\infty}
\frac{1}{\varrho(B^o_{T_i})}\int_{\tilde{A}^{C_0}_{T_i}}\varepsilon\omega(D_{s,T_i})e^{2\delta(s)}ds
=\varepsilon\lim_{T_i\rightarrow\infty} \frac{\varrho(\tilde{B}^{C_0}_{T_i})}{\varrho(B^o_{T_i})}=\varepsilon.
\end{eqnarray*}
This shows that $\eta(C)=0$. Hence, $\eta (Y)=0$.
\end{proof}\medbreak

Recall that by Lemma \ref{lem_eta}, $\eta$ is $U$-invariant.
By Ratner's measure classification \cite{r1}, an ergodic component of $\eta$ is
either $G$-invariant or supported on $Y\cup\{\infty\}$. By Lemmas
\ref{lem_inf1} and \ref{lem_Y1}, the set of ergodic components of the second type
has measure $0$. Therefore, $\eta$ is $G$-invariant, and $\eta=\nu$.
This finishes the proof of Theorem \ref{th_ergodic}.

\section*{ Acknowledgments }

The author wishes to thank his advisor V.~Bergelson for constant
support, encouragement, and help with this paper, N.~Shah and Barak~Weiss for insightful comments about the paper,
and H.~Furstenberg and G.~Margulis for interesting discussions.

\end{document}